\documentclass[11pt,twoside,a4paper]{article}

\usepackage{amsmath} \usepackage{amssymb} \usepackage{amsthm}
\usepackage{mathrsfs}

\usepackage[all]{xy} \usepackage{eucal} \SelectTips {cm}{}
\CompileMatrices

\usepackage{a4}
\newtheorem{thm}{Theorem}[section] \newtheorem{lem}[thm]{Lemma}
\newtheorem{prop}[thm]{Proposition} \newtheorem{cor}[thm]{Corollary}

\theoremstyle{definition}

\theoremstyle{remark} 

\newenvironment{absatz}[1][]
{\addtocounter{thm}{1} \noindent {\bf (\thesection.\arabic{thm}) #1}}{\rm}

\newenvironment{pro}[1][Proof]{{\it{#1:}} }{\hfill $\square$}
\newenvironment{pro*}[1][Proof]{{\it{#1:}} }{}

\numberwithin{equation}{section}

\newcommand{\mengensymb}[1]{\mathbb{#1}}
\newcommand{\A}{\mengensymb{A}} \newcommand{\C}{\mengensymb{C}}
 \newcommand{\N}{\mengensymb{N}}
\newcommand{\PP}{\mengensymb{P}} \newcommand{\Q}{\mengensymb{Q}}
 \newcommand{\Z}{\mengensymb{Z}}
 \newcommand{\LL}{\mengensymb{L}}

\newcommand{\OO}{\mathscr{O}}

\newcommand{\MM}{\mathscr{M}}
\DeclareMathOperator{\Moduli}{\mathscr{M}}

\newcommand{\Hom}{{\rm Hom}}
\newcommand{\Aut}{{\rm Aut}}
\newcommand{\Inn}{{\rm Inn}}
\newcommand{\Out}{{\rm Out}}
\newcommand{\Isom}{{\rm Isom}}
\newcommand{\surj}{\twoheadrightarrow}

\newcommand{\pll}{\pi^\LL}

\newcommand{\GMMgn}{\sideset{_G}{_{g,n}}\Moduli}

\DeclareMathOperator{\Spec}{{\rm Spec}}

\newcommand{\bru}[3]{\genfrac{}{}{#1}{}{#2}{#3}}
\newcommand{\ov}[1]{\mbox{${\overline{#1}}$}}

\newcommand{\ph}{\varphi}
\newcommand{\Ho}{{\rm H}}

\newcommand{\mgn}{\MM_{g,n}}
\newcommand{\mg}{\MM_{g}}
\newcommand{\mgo}{\MM_{g,0}}
\newcommand{\mzo}{\MM_{2,0}}
\newcommand{\mee}{\MM_{1,1}}

\newcommand{\modrei}{\MM_{0,3}}
\newcommand{\Jac}{{\rm Jac}}

\newcommand{\pr}{{\rm pr}}

\newcommand{\et}{\text{\rm \'et}}

\begin{document}

\title{$\phantom{.}$ \\[-6ex] \hrule $\phantom{a}$ \\[3ex] \bf A monodromy
  criterion for extending curves} \author{Jakob Stix\thanks{\sc Mathematisches
    Institut, Universit\"at Bonn, Beringstra\ss e 1, 53115 Bonn \newline
    \hspace*{0.45cm} E-mail address: {\tt stix@math.uni-bonn.de}\newline The
    author acknowledges the financial support provided through the European
    Community's Human Potential Program under contract HPRN-CT-2000-00114,
    GTEM.}}  \date{\today}

\maketitle

\begin{quotation} 
  \noindent \small {\bf Abstract} ---
  A family of proper smooth curves of genus $\geq 2$, parametrised by an open
  dense subset $U$ of a normal variety $S$, extends to $S$ if the natural map
  $\pi_1 U \to \pi_1 S$ is an isomorphism. The criterion of this note is
  actually more precise. 
\end{quotation}


\section{The monodromy criterion}

\begin{absatz}
  All schemes are locally noetherian. For a scheme $S$ we denote by $\mgn(S)$
  the groupoid of smooth, proper $S$-curves with geometrically connected
  fibres of genus $g$ and $n$ disjoint ordered sections, see \cite{Knudsen}.
  If $2-2g-n$ is negative, $\mgn$ is a smooth Deligne-Mumford stack over
  $\Spec(\Z)$. The main result of this note is the following theorem.
\end{absatz}

\begin{thm}[monodromy criterion I] \label{t}
  Let $2-2g-n$ be negative.  For a dense open subscheme $U$ of a normal,
  excellent scheme $S$ the following holds.
\begin{itemize}
\item[(1)] The restriction functor $\mgn(S) \to \mgn(U)$ is fully faithful.
\item[(2)] A $U$-curve $C \in \mgn(U)$ extends to an $S$-curve in $\mgn(S)$ if
  and only if for all $N \in \N$ and geometric points $u \in U \otimes
  \Z[\bru{0.5pt}{1}{N}]$ the following commutative diagram exists, where the
  solid arrows are induced by the natural maps $U \to \mgn$ and $U \subset S$:
  $$
  \xymatrix{{\pi_1(U \otimes \Z[\bru{0.5pt}{1}{N}],u)} \ar[rr] \ar[dr] &&
    {\pi_1(\mgn \otimes \Z[\bru{0.5pt}{1}{N}],C_u)} \\
    & {\pi_1(S \otimes \Z[\bru{0.5pt}{1}{N}],u)} \ar@{.>}[ur]& }$$
\end{itemize}
\end{thm}
\begin{absatz}[Remarks.] (a) The base point $C_u$ of $\mgn$ is the
  isomorphy class of the fibre $C_u$ of $C$ over $u$. The fundamental group of
  the stack $\mgn$ is the pro-finite group classifying finite \'etale covers
  of $\mgn$.
  
  (b) We may and will assume that $S$ and consequently $U$ are connected.
  Furthermore we may and will ignore base points. The map $\pi_1(U \otimes
  \Z[\bru{0.5pt}{1}{N}]) \to \pi_1(S \otimes \Z[\bru{0.5pt}{1}{N}])$ is
  surjective for all $N$, such that $U \otimes \Z[\bru{0.5pt}{1}{N}]$ is not
  empty.
  
  (c) We will discuss later that Theorem \ref{t} follows from the work of
  Moret-Bailly and Oda-Tamagawa in case $S$ is regular. Thus the new
  contribution of the present note consists in the case of a singular, but
  normal base $S$.
\end{absatz}

\begin{absatz}
  For $C,C' \in \mgn(S)$ a modification to the pointed case of \cite{DM} Thm
  1.11 shows that the $S$-scheme $\Isom(C/S,C'/S)$ is finite (unramified) over
  $S$, i.e., the stack $\mgn$ is a separated Deligne--Mumford stack, see
  \cite{Knudsen}.  Hence a $U$-valued point for a dense open $U$ of a normal
  scheme $S$ extends uniquely to an $S$-valued point. This proves part (1) of
  Theorem \ref{t}.
 
  We may and will replace $S$ by $S \otimes \Z[\bru{0.5pt}{1}{N}]$ for various
  $N \in \N$ and glue afterwards using Theorem \ref{t} (1). In particular we
  may assume that some prime $\ell \in \N$ is invertible in $S$.
\end{absatz}
\begin{cor}
  Let $2-2g-n$ be negative and let $k$ be a field.  Let $U$ be a dense open
  subscheme of a normal, connected, excellent scheme $S$ over $\Spec(k)$.
  
  Then a $U$-curve $C \in \mgn(U)$ extends to an $S$-curve in $\mgn(S)$ if and
  only if the following commutative diagram exists, where the solid arrows are
  induced by the natural maps $U \to \mgn \otimes k$ and $U \subset S$:
  $$
  \xymatrix{{\pi_1U} \ar[rr] \ar[dr] &&
    {\pi_1(\mgn \otimes k)} \\
    & {\pi_1S} \ar@{.>}[ur]& }$$
  \vspace{-1cm}
  
  \ \hfill $\square$
\end{cor}

\begin{absatz}[Remark.]
  The motivation behind Theorem \ref{t} stems from Grothendieck's idea that
  $\mgn$ behaves `anabelian', meaning that the pro-finite \'etale fundamental
  group should control the geometry of $\mgn$. In the analytic context, the
  orbifold space $\mgn(\C)$ is a ${\rm K}(\pi,1)$ space for the Teichm\"uller
  mapping class group. Hence for varieties $S$ over $\C$ our Theorem \ref{t}
  is trivial up to homotopy by topological methods.  But still one needs to
  prove that some homotopy classes of continuous maps from $S(\C)$ to
  $\mgn(\C)$ contain a unique algebraic map extending a given map on the dense
  open $U(\C)$.
  
  Nevertheless, the difference between the analytic and algebraic situation
  lies in the pro-finite completion of the fundamental group. The good ${\rm
    K}(\pi,1)$ properties only prevail in the pro-finite setting if the
  mapping class group is a good group in the sense of Serre. And this is not
  known.
\end{absatz}


\section{The outer monodromy representation} \label{secOR}

\begin{absatz}
  In what follows, $\LL$ is a set of prime numbers invertible on $S$.  Let
  $C/S$ be an $S$-curve in $\mgn(S)$. Let $C^0$ be the notation for the open
  complement of its sections. Then $C^0_s$ is the fibre of $C^0/S$ over the
  geometric point $s \in S$ and $\pll_{s,c} = \pi_1^\LL (C^0_{s,c})$ the
  pro-$\LL$ completion of its fundamental group, i.e., the maximal continuous
  pro-finite quotient that is a pro-$\LL$ group (a pro-finite limit of groups
  of order a product of powers of primes from $\LL$). The group $\pll_{s,c}$
  is noncanonically isomorphic to the pro-$\LL$ completion $\pll$ of the
  topological fundamental group of an oriented topological surface of genus
  $g$ with $n$ cusps.
\end{absatz}
\begin{lem} \label{lemPL}
  (1) The group $\pll$ is finitely generated as a pro-finite group. The groups
  $\Aut(\pll)$ and $\Out(\pll)$ are canonically pro-finite groups.  (2) If
  $2-2g-n$ is negative, then the center of $\pll$ is trivial.
\end{lem}
\begin{pro}
  (1) A finitely generated pro-finite group is a pro-finite limit of
  characteristic finite quotients. (2) See \cite{Anderson} Prop 8 and Prop 18.
\end{pro}

\begin{absatz}
  Let $K$ denote the kernel of the map $\pi_1(C^0,c) \to \pi_1(S,s)$ and $N$
  be its smallest normal closed subgroup such that $K/N$ is a pro-$\LL$ group.
  Then $N$ is characteristic in $K$ and thus also a normal subgroup of
  $\pi_1(C^0,c)$. We set $\pi_1^{(\LL)}(C^0,c)$ for the quotient
  $\pi_1(C^0,c)/N$. The sequence
\begin{equation*} 
\pll_{s,c} \to \pi_1^{(\LL)}(C^0,c) \to \pi_1(S,s) \to 1 \tag{$\ast$}
\end{equation*}
is exact by \cite{SGA1} XIII 4.1/4.4. Moreover, if $C^0/S$ admits a section
$\sigma$, then the sequence ($\ast$) is even split exact by \cite{SGA1} XIII
4.3/4.4, the splitting being induced by the section.
\begin{equation*} 
\xymatrix@1{1 \ar[r] & \pll_{s,\sigma(s)}  \ar[r] &
 \pi_1^{(\LL)}(C^0,\sigma(s)) \ar[r]  & \pi_1(S,s) \ar@{.>}@/_2ex/[l]_\sigma
 \ar[r] & 1 } \tag{$\ast\ast$}
\end{equation*}
Up to a choice of an isomorphism $\pll \cong \pll_{s,\sigma(s)}$ the
corresponding outer pro-$\LL$ representation
$$
\rho_{C/S} : \pi_1(S,s) \to \Out\big(\pll\big)$$
induced by conjugation is
independent of the section and obeys functoriality with respect to maps $f:S'
\to S$ and base change $f^\ast C = C \times_S S'$, i.e., $\rho_{f^\ast C/S'} =
\rho_{C/S} \circ \pi_1f$, when regarded as an outer homomorphims (up to
composition with an inner automorphism).
\end{absatz} 

\begin{absatz}
  Changing the isomorphism $\pll \cong \pll_{s,\sigma(s)}$ composes
  $\rho_{C/S}$ by an inner automorphism of $\Out(\pll)$, thus leaving the
  class as an outer homomorphism unchanged. In the sequel outer
  representations will be regarded as outer homomorphisms. Diagrams of outer
  homomorphisms commute if they `commute up to inner automorphisms' for
  representatives.
\end{absatz}

\begin{absatz}
  In fact, the canonical outer pro-$\LL$ representation does also exist in
  absence of a section by the following proposition; see \cite{MT} for a
  characteristic $0$ analogue.
\end{absatz}
\begin{prop} \label{propRho}
  Let $S$ be a connected scheme and let $C/S$ be an $S$-curve in $\mgn(S)$.
  Let $\LL$ be a set of prime numbers which are invertible on $S$.
\begin{itemize}
\item[(1)] There is a unique outer pro-$\LL$ representation $\rho :\pi_1S \to
  \Out\big(\pll\big)$, such that for all $f: S' \to S$ where $f^\ast C/S'$
  admits a section the following diagram commutes.
  $$
  \xymatrix{{\pi_1S'} \ar[rr]^{\rho_{f^\ast C/S'}} \ar[dr]_{\pi_1f} &&
    {\Out\big(\pll\big)} \\
    & {\pi_1S} \ar[ur]_{\rho}& }$$
\item[(2)] For all $f:S' \to S$ with $S',S$ connected the outer pro-$\LL$
  representations $\rho$ (resp.  $\rho'$) from (1) applied to $C/S$ (resp.
  $f^\ast C/S'$) satisfy $\rho' = \rho \circ \pi_1f$.
\end{itemize}
\end{prop} 
\begin{pro}
  (2) follows from the uniqueness in (1) which follows from ($\ast$) and the
  diagonal section of $C^0 \times_S C^0$. As any base change which allows a
  section of $C^0/S$ factors through $\pr: C^0 \to S$, it suffices for (1) to
  prove that $\rho'=\rho_{\pr^\ast C/C^0}$ factors through $\pi_1C^0 \surj
  \pi_1S$.
  
  We follow the construction in \cite{Pikaart-deJong} \S2.1. For simplicity
  let us assume for a moment that $C^0/S$ has a section $\sigma$. The outer
  representation $\rho_{C/S}$ lifts via $\sigma$ to an actual continuous
  morphism $\pi_1(S,s) \to \Aut(\pll_{s,\sigma(s)})$ that describes a
  pro-object of locally constant sheaves of groups on $S_\et$, denoted
  $\pi_1^\LL(C^0/S,\sigma)$, with fibre $\pll_{s,\sigma(s)}$ in $s$.
  
  The sheaf $\Isom^{\rm out}\big(\pll,\pi_1^\LL(C^0/S,\sigma)\big)$ of outer
  isomorphisms of the constant pro-object $\pll$ with
  $\pi_1^\LL(C^0/S,\sigma)$ forms itself a pro-object of locally constant
  sheaves of sets: more precisely an $\Out(\pll)$-torsor by composition (from
  the right) with fibre $\Isom^{\rm out}\big(\pll,\pll_{s,\sigma(s)}\big)$ in
  $s$. The choice $t$ of such an isomorphism $\pll \cong \pll_{s,\sigma(s)}$
  allows us to describe the torsor as an actual homomorphism $\pi_1(S,s) \to
  \Out(\pll)$ which is nothing but $\rho_{C/S}$ again.
  
  For different sections $\sigma$ and $\sigma'$ there are isomorphisms
  $\alpha$ an $\alpha^\LL$ that induce the following isomorphism of extensions
  that is not necessarily compatible with the sections.
  $$
  \xymatrix{ 1 \ar[r] & \pll_{s,\sigma(c)} \ar[r] \ar[d]^{\alpha^\LL} &
    \pi_1^{(\LL)}(C^0,\sigma(c)) \ar[r] \ar[d]^\alpha & \pi_1(S,s) \ar@{=}[d]
    \ar[r] \ar@{.>}@/_2ex/[l]_(0.4)\sigma &
    1\\
    1 \ar[r] & \pll_{s,\sigma'(c)} \ar[r] & \pi_1^{(\LL)}(C^0,\sigma'(c))
    \ar[r] & \pi_1(S,s) \ar@{.>}@/_2ex/[l]_(0.4){\sigma'} \ar[r] & 1}
  $$
  Nevertheless, $\alpha$ and $\alpha^\LL$ are unique up to composition with
  an inner automorphism from an element of $\pll_{s,\sigma'(s)}$. In
  particular, the composed isomorphism
  $$\ph = t'^{-1} \alpha^\LL t : \pll \cong \pll_{s,\sigma(s)} \cong
  \pll_{s,\sigma'(s)} \cong \pll$$
  is unique in $\Out(\pll)$ and yields a
  canonical isomorphism
  $$\circ \ph : \Isom^{\rm out}\big(\pll,\pi_1^\LL(C^0/S,\sigma)\big)
  \xrightarrow{\sim} \Isom^{\rm out}\big(\pll,\pi_1^\LL(C^0/S,\sigma')\big)$$
  of $\Out(\pll)$-torsors.
  
  Now we return to the general case $C/S$ where $C^0/S$ need not have a
  section. We show that $\rho_{\pr^\ast C/S} : \pi_1(C^0) \to \Out(\pll)$,
  which is given by the diagonal section $\Delta$, factors through $\pi_1C^0
  \surj \pi_1S$. For this it is necessary and sufficient to descend the
  $\Out(\pll)$-torsor $\Isom^{\rm out}\big(\pll,\pi_1^\LL((\pr^\ast
  C)^0/C^0,\Delta)\big)$ from $C^0$ to $S$.  A descent datum is given by
  independence of the section up to canonical isomorphism. Finite \'etale
  torsors obey effective fpqc-descent, here along $C^0 \to S$, thus we obtain
  $\rho: \pi_1(S,s) \to \Out(\pll)$ as required.
  
  Changing the basepoint $s$ or changing the ismorphism $\pll \cong
  \pll_{s,\sigma(s)}$ affects the outer representation $\rho$ by composing
  with an inner automorphism of $\Out(\pll)$. This concludes the proof of the
  proposition.
\end{pro}

\begin{prop} \label{propEXT}
  Let $2-2g-n$ be negative. Let $S$ be a connected scheme and let $C/S$ be an
  $S$-curve in $\mgn(S)$. Let $\LL$ be a set of prime numbers which are
  invertible on $S$.  The following diagram has exact rows.
  $$\xymatrix{ 1 \ar[r] & \pll_{s,c} \ar[r] \ar[d]^{\cong} &
    \pi_1^{(\LL)}(C^0,c) \ar[r] \ar[d]^{\tilde{\rho}} & \pi_1(S,s) 
    \ar[d]^{\rho} \ar[r] & 1 \\
    1 \ar[r] & \Inn(\pll) \ar[r] & \Aut(\pll) \ar[r] & \Out(\pll) \ar[r] &
    1}$$
  The extension of groups in the upper row is induced via $\rho$ by the
  extension of the lower row.
\end{prop}
\begin{pro}
  It suffices to prove that $\pll_{s,c} \to \pi_1^{(\LL)}(C^0,c)$ is
  injective. For then the outer representation $\rho$ of Proposition
  \ref{propRho} equals the conjugation action by lifts under a chosen
  identification $\pll \cong \pll_{s,c}$. In particular the map $\tilde{\rho}$
  is induced by conjugation.
 
  Let $K$ be as above the kernel of $\pi_1(C^0,c) \to \pi_1(S,s)$.  Because
  $\rho_{\pr^\ast C/C^0}$ factors through $\pi_1(C^0,c) \to \pi_1(S,s)$ by
  Proposition \ref{propRho}, the lift $\tilde{\rho} \circ \Delta: \pi_1(C^0,c)
  \to \Aut(\pll_{s,c})$ invoking the diagonal section $\Delta$ of $(\pr^\ast
  C)^0/C^0$ as above and the conjugation action $\tilde{\rho}$ maps $K$ to the
  group of inner automorphisms $\Inn(\pll_{s,c})$ which is a pro-$\LL$ group.
  In particular $\tilde{\rho} \circ \Delta$ factors through a map
  $\pi_1^{(\LL)}(C^0,c) \to \Aut(\pll_{s,c})$. The restriction of
  $\tilde{\rho}\circ\Delta$ to $\pll_{s,c}$ coincides with the canonical map
  $\pll_{s,c} \cong \Inn(\pll_{s,c})$ being induced by the diagonal section of
  $$
  1 \to \pll_{s,c} \to \pll_{s,c} \times \pll_{s,c} \to \pll_{s,c} \to 1\ .
  $$
  Here we used a pro-$\LL$ product formula for $\pi_1$ which is valid as
  all primes from $\LL$ are invertible on $S$. The vanishing of the center of
  $\pll$, see Lemma \ref{lemPL} (2), shows that $\tilde{\rho} \circ \Delta$ is
  injective on $\pll_{s,c}$ and thus also the injectivity of the map
  $\pll_{s,c} \to \pi_1^{\LL}(C^0,c)$.
\end{pro}

\begin{absatz}[Remarks.]
  (a) Let $\LL$ be a set of prime numbers. We set $\Z[\bru{0.5pt}{1}{\LL}]$
  for the ring of rational numbers whose denominator only contains primes from
  $\LL$.  The construction of the $\Out(\pll)$-torsor in Proposition
  \ref{propRho} may be applied to the universal curve above $\mgn \otimes
  \Z[\bru{0.5pt}{1}{\LL}]$. We obtain the universal outer pro-$\LL$
  representation
  $$
  \rho^{\rm univ} : \xymatrix@1{\pi_1 \big(\mgn \otimes
    \Z[\bru{0.5pt}{1}{\LL}]\big) \ar[r] & \Out(\pll)}.$$
  
  (b) Let $G$ be a finite group. We may fuse the $\Out(\pll)$-torsor with the
  $\Out(\pll)$-set $\Hom_{\rm surj}^{\rm out}(\pll,G)$ of surjective, outer
  homomorphisms $\pll \surj G$. We obtain a finite \'etale cover $\GMMgn \to
  \mgn\otimes \Z[\bru{0.5pt}{1}{\#G}]$ of curves that are equipped with a
  nonabelian Teichm\"uller level structure of level $G$, see
  \cite{Pikaart-deJong} \S 2.3 and \cite{BoggiPikaart}. The stack $\GMMgn$ is
  an open substack of the normalisation $\ov{\GMMgn}$ of the stack
  $\ov{\MM}_{g,n}$ of stable pointed curves in the function field of $\GMMgn$.
  
  For $d$ a product of powers of primes from $\LL$ and $G= \pll[d] :=
  \pll/[\pll,\pll](\pll)^d$ we denote the space $\GMMgn$ of curves with
  abelian level structures of level $d$ by $\mgn[d]$.
  
  If $G$ surjects to $\pll[\ell^r]$ for some $\ell \in \LL$ and $\ell^r\geq
  3$, then $\ov{\GMMgn}$ is actually a scheme, see \cite{Pikaart-deJong} Prop
  2.3.4, which is projective over $\Spec(\Z[\bru{0.5pt}{1}{\#G}])$.
  
  (c) Note that the assumption $2-2g-n<0$ in Proposition \ref{propEXT} is
  essential. A counterexample for the case $\MM_{0,2}$ is as follows.  The
  complement of the sections of the $\PP^1$-bundle $\PP(\OO(-1)\oplus \OO)$
  over $\PP_\C^1$ equipped with the two sections $0$ and $\infty$ is
  isomorphic to $\A_\C^2-\{0\}$. Hence it is simply connected, whereas a fibre
  has fundamental group isomorphic to $\widehat{\Z}$.
\end{absatz}

\begin{absatz}
  We can formulate a second monodromy criterion for extending curves in terms
  of outer representations as follows.
\end{absatz}
\begin{thm}[monodromy criterion II] \label{to}
  Let $2-2g-n$ be negative. Let $U$ be a dense open subscheme of a normal,
  connected, excellent scheme $S$.
  
  Let $C/U$ be a $U$-curve in $\mgn(U)$ and let $\pll$ be as above the
  pro-$\LL$ completion of the complement of the sections of $C/U$ in a
  geometric fibre. Then $C/U$ extends to an $S$-curve in $\mgn(S)$ if and only
  if the following commutative diagram, where the solid arrows are induced by
  the natural maps $U \subset S$ and the outer pro-$\LL$ representation $\rho$
  associated to $C\otimes \Z[\bru{0.5pt}{1}{\LL}]/U\otimes
  \Z[\bru{0.5pt}{1}{\LL}]$,
  $$
  \xymatrix{{\pi_1\big(U\otimes \Z[\bru{0.5pt}{1}{\LL}]\big)} \ar[rr]^\rho
    \ar[dr] &&
    {\Out\big(\pll)} \\
    & {\pi_1\big(S\otimes\Z[\bru{0.5pt}{1}{\LL}]\big)} \ar@{.>}[ur] & }$$
  exists for one of the following collections of sets of prime numbers $\LL$
  and additional conditions on $S$:
\begin{itemize}
\item[(A)] $\LL=\{\ell\}$ for some prime number $\ell$ and $S$ is of
  characteristic $0$, or
\item[(B)] $\LL=\{\ell_1, \ell_2\}$ for all pairs of sufficiently large prime
  numbers $\ell_1,\ell_2$ and no additional conditions on $S$.
\end{itemize} 
\end{thm}
\begin{absatz}[Remarks.] \label{rmkOT}
  (a) It follows from the universal outer pro-$\LL$ representation $\rho^{\rm
    univ}$ and functoriality that Thereom \ref{to} is stronger than Theorem
  \ref{t}.
  
  (b) It is plausible though, that Theorem \ref{to} also holds under the
  following condition on the collection of $\LL$'s:
\begin{itemize}
\item[(C)] $\LL=\{\ell\}$ for some prime number $\ell$ that is invertible on
  $S$.
\end{itemize}
Indeed, assuption (B) which involves two prime numbers only comes into play to
apply results from \cite{Tama2} in the positive characteristic situation of
Theorem \ref{tc}, where alternative proofs are not out of sight, see the
remarks there.
  
(c) In case $S$ is the spectrum of a henselian discrete valuation ring,
Theorem \ref{to}, under condition (C) above, is nothing but the criterion for
good reduction of Oda--Tamagawa \cite{Tama1} Thm 5.3.  In particular, curves
satisfying either monodromy criterion above extend into codimension $1$ points
of the base. It remains to deal with $U \subset S$ with boundary of
codimension $\geq 2$.

(d) In \cite{MoretBailly} Moret-Bailly proves a purity theorem for relative
curves over regular bases $S$. Namely, for $U$ open dense in a regular $S$
with $S-U$ of codimension at least $2$, the restriction functor $\mgn(S) \to
\mgn(U)$ is an equivalence of categories. Moret-Bailly's theorem together with
Oda-Tamagawa's criterion as mentioned in (c) imply our Theorem \ref{to}, if
the base $S$ is regular along the boundary $S\setminus U$. In fact, by
Zariski--Nagata purity for the branch locus, Moret-Bailly's theorem is
equivalent to our monodromy criteria in the case of a regular base and
codimension of the boundary at least $2$.

The new contribution of the present note thus only consists in the case where
$S$ is not regular along the boundary.

In the beginning, the proof in \cite{MoretBailly} and our treatment of the
monodromy criterion Theorem \ref{to} follow the same strategy while the
methods to enforce isotriviality along subvarieties, that one needs to
contract, are essentially different. More precisely, \cite{MoretBailly} uses
first a deformation argument to reduce to bases of dimension $2$ to the effect
that by induction and the theory of regular surfaces it suffices to contract
rational lines. Then \cite{MoretBailly} uses results on pencils of stable
curves over $\PP^1$.

The deformation argument is not applicable in our case where the base is only
assumed to be normal.  Instead, the present note exploits a recent result of
Tamagawa \cite{Tama2} on the specialisation homomorphism for the fundamental
group of curves in positive characteristic, see also \cite{Saidi}.
  
(e) The result of \cite{MoretBailly} has been extended in various ways to
stable or log smooth curves by de~Jong--Oort \cite{dJongOort}, Mochizuki
\cite{Mochiz} and T.~Saito \cite{Saito2}. The author believes that the
monodromy criteria of this note extend to the logarithmic situation, cf.
\cite{stix} Thm 1.2.
  
(f) Theorem \ref{ext} at the end of this note contains an extended monodromy
criterion relative certain maps of finite type.
\end{absatz}  


\section{Extension after modification of the base}

\begin{absatz}
  The rest of this note is devoted to the proof of Theorem \ref{to} and thus
  of part (2) of Theorem \ref{t}. We only need to show that a curve $C/U$
  extends if the condition on the fundamental groups holds.
\end{absatz}

\begin{absatz}
  We can use part (1) of Theorem \ref{t} and \'etale descent for curves to
  work locally on $S$ in the \'etale topology. The descent is effective by the
  presence of a canonical, relatively ample line bundle: the bundle of
  relative differentials twisted by the locus of the sections.
  
  For example, we may replace $S$ by a finite \'etale cover $S'$ obtained by
  pullback of a finite \'etale cover $\Moduli'$ of $\mgn \otimes
  \Z[\bru{0.5pt}{1}{\ell}]$ that is induced via the universal outer pro-$\ell$
  representation by a finite quotient of $\Out\big(\pi^\ell\big)$. Indeed,
  under the map $U \to \mgn$ the cover $\Moduli'$ pulls back to a finite
  \'etale cover $U'$ of $U$ that extends to $S'$ by the condition on the
  $\pi_1$'s.
  
  We apply the preceeding construction to the cover $\mgn[\ell^r]$ of curves
  with abelian level structure of level $\ell^r \geq 3$ for some prime $\ell
  \in \LL$.
\end{absatz}

\begin{absatz}
  The first step will be to prove that if a curve in $\mgn[\ell^r](U)$
  satisfies the monodromy criterion then it extends as a curve with abelian
  level structure to a proper modification $S' \to S$ which is an isomorphism
  over $U$.  This part of the theorem holds unconditionally on $\ell \in
  \OO_S^\ast$.
\end{absatz}

\begin{thm}[extension after modification] \label{tm}
  Let $U$ be an open dense subscheme in a normal, connected, excellent scheme
  $S$, and let $\ell$ be a prime number invertible on $S$. Let $C/U$ be a
  $U$-curve in $\mgn[\ell^r](U)$ endowed with an abelian level structure of
  level $\ell^r$ for some $\ell^r \geq 3$.
  
  If the associated outer pro-$\ell$ representation factors through $\pi_1U
  \to \pi_1S$, then there exists a proper birational map $\sigma : S' \to S$
  which is an isomorphism above $U$, such that
\begin{itemize}
\item[(i)] $S'$ is normal,
\item[(ii)] $\sigma_\ast \OO_{S'} = \OO_S$,
\item[(iii)] the $U$-curve $C/U$ extends to an $S'$-curve with level structure
  in $\mgn[\ell^r](S')$ such that the corresponding outer pro-$\ell$
  representation factors through $\pi_1\sigma : \pi_1 S' \to \pi_1S$.
\end{itemize}
Moreover, if also the outer pro-$\LL$ representation of $C/U$ factors through
$\pi_1U \to \pi_1S$ for some set of primes $\LL$ invertible on $S$, then the
correponding outer pro-$\LL$ representation of the extension over $S'$ factors
through $\pi_1\sigma : \pi_1 S' \to \pi_1S$.
\end{thm}
\begin{pro}
  The curve $C/U$ is represented by a map $f_U:U \to \mgn[\ell^r]$ and we need
  to find a modification $\sigma : S' \to S$ such that $f_U$ extends to a map
  $f':S' \to \mgn[\ell^r]$. Let $S'$ be the normalisation of the closure in $S
  \times \ov{\Moduli}_{g,n}[\ell^r]$ of the graph of $f_U$.  The first
  projection $\sigma : S' \to S$ satisfies the requirements.  To see this, we
  argue that the second projection $f' : S' \to \ov{\Moduli}_{g,n}[\ell^r]$
  actually has image contained in $\mgn[\ell^r]$. The other requirements for
  $\sigma$ are immediate.
  
  Assume on the contrary that $s' \in S'$ is mapped to $f'(s')$ in the
  boundary. We can find a strict henselian, discrete valuation ring $R$ and
  $\Spec(R) \to S'$ with the closed point mapping to $s'$ and the generic
  point $\eta$ mapping into $U \subset S'$. The composition with $f'$ yields a
  map $\Spec(R) \to \ov{\Moduli}_{g,n}[\ell^r]$ such that the corresponding
  curve over $\Spec(R)$ has bad reduction.  The assumption on the outer
  pro-$\ell$ representation implies that the corresponding map $\pi_1\eta \to
  \pi_1U \to \Out\big(\pi^\ell\big)$ is trivial, This contradicts the
  Oda--Tamagawa criterion for good reduction, \cite{Tama1} Thm 5.3, that was
  already mentioned above in remark (2.11 (c)).
\end{pro}

\begin{absatz}[Remark.]
  For $g\geq 2$ and $n=0$ it follows from \cite{Pikaart-deJong} Thm 3.1.3 that
  $\mgn \otimes \Z[\bru{0.5pt}{1}{\ell}]$ is universally ramified in
  $\ov{\Moduli}_{g,n} \otimes \Z[\bru{0.5pt}{1}{\ell}]$ with respect to
  $\ell$, see also Brylinski, Looijenga from the references of loc.\ cit.
  Conjecturally $\mgn \otimes \Z[\bru{0.5pt}{1}{\ell}]$ has universal
  ramification in $\ov{\Moduli}_{g,n} \otimes \Z[\bru{0.5pt}{1}{\ell}]$ with
  respect to $\ell$ for all $(g,n)$ with $2-2g-n$ negative. There is a proof
  of Theorem \ref{tm} based on universal ramification, see \cite{stix2} \S2.
  This would lead to a geometric proof based on the moduli space of the
  Oda--Tamagawa criterion for good reduction as a special case of Theorem
  \ref{to}.
\end{absatz}


\section{Constant maps} \label{secCM}

\begin{absatz}
  To complete the proof of Theorem \ref{to} it remains to prove that the map
  $f'$ from the proof of Theorem \ref{tm} factors through $\sigma : S' \to S$
  under the monodromy conditions of Theorem \ref{to}. More precisely, the
  latter imply as indicated in Theorem \ref{tm}, that the respective outer
  pro-$\LL$ representations $\pi_1S' \to \Out(\pll)$ factor through
  $\pi_1\sigma : \pi_1S' \to \pi_1S$.
\end{absatz}

\begin{absatz}
  It suffices that $f'$ factorises set-theoretically as $f'=f \circ \sigma$.
  Indeed, $S$ carries the quotient topology under $\sigma$, hence $f$ is
  automatically continuous. Then the property $\sigma_\ast\OO_{S'}=\OO_S$ is
  responsible for the enrichment of $f$ to a map of (locally) ringed spaces.
\end{absatz}

\begin{absatz}
  For $g \leq 2$ there is nothing to prove. Indeed, $\mgn[\ell^r]$ is
  consecutively fibred in affine curves over $\mzo[\ell^r]$, $\mee[\ell^r]$ or
  $\modrei[\ell^r]$ respectively. These latter moduli spaces are affine,
  whereas the fibres of $\sigma$ are proper. In particular Theorem \ref{to}
  holds even under hypothesis (C) on $\ell$ and $S$ in these cases.
\end{absatz}

\begin{absatz}
  For $g > 2$ and a geometric fibre $F$ of $\sigma$, the outer pro-$\LL$
  representations $\rho: \pi_1F \to \Out\big(\pll\big)$ of the restriction of
  the $S'$-curve to $F$ factors through $\pi_1S$ for the various sets $\LL$ in
  question. Hence $\rho$ is trivial.
  
  The following Theorem \ref{tc} accomplishes the last step in the proof of
  Theorem \ref{to} and thus Theorem \ref{t}.
\end{absatz}

\begin{thm} \label{tc}
  Let $2-2g-n$ be negative.  Let $F$ be a reduced, connected variety over an
  algebraically closed field $k$. Let $\ph : F \to \mgn$ be a map such that
  the outer pro-$\LL$ representation $\rho : \pi_1F \to \Out\big(\pll\big)$ is
  the trivial homomorphism for one of the following collections of sets of
  prime numbers $\LL$ and additional conditions on $k$:
\begin{itemize}
\item[(A)] $\LL=\{\ell\}$ for some prime number $\ell$ and $k$ is of
  characteristic $0$, or
\item[(B)] $\LL=\{\ell_1, \ell_2\}$ for all pairs of sufficiently large prime
  numbers $\ell_1,\ell_2$ invertible in $k$ and $k$ is of positive
  characteristic.
\end{itemize} 
Then $\ph$ is constant in the sense that the correspondig $F$-curve $C/F \in
\mgn(F)$ comes by base extension from an object in $\mgn\big(\Spec(k)\big)$.
\end{thm}

\begin{absatz}[Remark.] The characteristic $0$ part of the above theorem
  follows also from the work of Mochizuki on pro-$p$ anabelian geometry, see
  \cite{Mz}.
\end{absatz}

\begin{pro}
  Let $\ell$ be a prime in $\LL$.  By assumption, the map $\ph$ lifts to a map
  $\tilde{\ph}:F \to \mgn[\ell^r]$. As for $\ell^r\geq 3$ the target is now a
  scheme, it suffices to prove that $\tilde{\ph}$ is a constant map of sets in
  this case.
  
  By covering $F$ with images of curves we may assume that $F$ is a smooth
  curve. Now we apply Theorem \ref{tm} to $\tilde{\ph} : F \to \mgn[\ell^r]$
  and deduce that we may extend the curve to the smooth compactification of
  $F$. Still, the outer pro-$\LL$ representation is trivial for the respective
  $\LL$.  Hence we may assume that $F$ is even a proper smooth curve over $k$.
  
  Next we observe that $\MM_{g,n+1}[\ell^r] \to \mgn[\ell^r]$ is fibred in
  affine curves for $n \geq 1$. By induction on $n$ we are reduced to the
  cases $n \leq 1$. The last step $\MM_{g,1}[\ell^r] \to \mgo[\ell^r]$ is
  fibred in smooth proper curves of fixed genus $g$. Let us assume that we
  have successfully dealt with the case $n=0$. Then any non constant lift
  corresponds to a dominant map from $F$ to a curve of genus $g$ that moreover
  lifts by the assumption on $\rho$ to any finite $\ell$-primary Galois cover.
  As the genus of such covers tends to infinity (here $g \geq 2$) this is
  impossible.
  
  It remains to deal with the cases $g \geq 3$ and $n=0$, because the moduli
  spaces are affine in the other minimal cases $(2,0)$, $(1,1)$ and $(0,3)$.
  
  We fix geometric points $x \in F$ and $c \in C_x$, the fibre of $C/F$ above
  $x$. The vanishing of the outer pro-$\LL$ representation and Proposition
  \ref{propEXT} yield a canonical isomorphism
  $$\pi_1^{(\LL)}(C,c) = \pll_{x,c} \times \pi_1(F,x).$$
  The finite \'etale
  cover $C_{x,\Phi}$ of $C_x$ associated to a finite set $\Phi$ with
  continuous $\pll_{x,c}$-action is the fibre above $x$ of a finite \'etale
  cover $C_\Phi \to C$ associated to the action of $\pi_1^{(\LL)}(C,c)$ on
  $\Phi$ via projection to $\pll_{x,c}$.  If $C_{x,\Phi}$ is connected, then
  $C_\Phi$ is a family of geometrically connected curves over $F$. By
  construction we then have canonically $\pi_1^{(\LL)}C_\Phi = \pll_\Phi
  \times \pi_1F$ with the group $\pll_\Phi = \pi_1^\LL(C_{x,\Phi}) \subset
  \pll_{x,c}$.  It follows that any $C_\Phi/F$ as above with geometrically
  connected fibres has constant abelian $\ell^n$ level structure for all $n
  \in \N$ and $\ell \in \LL$ and therefore obeys the following theorem.

\begin{thm} \label{tg}
  Let $B$ be a normal variety over an algebraically closed field $k$.  Let
  $\ell$ be invertible in $k$. Let $C/B$ be a curve with constant abelian
  $\ell^n$ level structure for all $n \in \N$. Then the following holds.
\begin{itemize}
\item[(i)] The relative jacobian $\Jac_{C/B}$ is radicially isogenous to a
  constant abelian variety.
\item[(ii)] If $k$ has characteristic $p > 0$, then the $p$-rank of the fibres
  $C_b$ of $C/B$ is a constant function of $b \in B$.
\end{itemize}
\end{thm}

\begin{pro}
  Part (i) follows from \cite{Groth} Prop.\ 4.4 (or \cite{Oort} Thm 2.1), the
  essential ingredients are the theory of the $k$-trace of abelian varieties
  and the relative Mordell--Weil Theorem of Lang--N\'eron.  For (ii) just
  notice that isogenous abelian varieties have the same $p$-rank.
\end{pro}

\begin{absatz}
  We come back to the proof of Theorem \ref{tc}.  If $k$ is of characteristic
  $0$ then by (i) of Theorem \ref{tg} the relative jacobian together with its
  canonical principal polarization and abelian $\ell^n$ level structure is
  constant. By Torelli's theorem, see \cite{OS}, the same holds for the curve
  $C/F$ and we are done. It remains to deal with positive characteristic.
\end{absatz}

\begin{absatz}
  Let $k$ have positive characteristic $p$. We deduce from Theorem \ref{tg}
  that all families $C_\Phi/F$ as above have constant $p$-rank.
\end{absatz}

\begin{absatz}
  Following ideas of Raynaud (for new ordinarity), in \cite{Tama2} Theorem
  0.7, Proposition 0.8, Corollary 5.3, and Theorem 0.5 Tamagawa proved the
  existence of prime to $p$ covers $C''_x \to C'_x \to C_x$ of $C_x = C
  \times_F x$, such that $C''_x/C'_x$ is
\begin{itemize}
\item[(i)] prime cyclic Galois,
\item[(ii)] new ordinary: the growth in $p$-rank equals the growth in genus,
  and
\item[(iii)] `new Torelli': the infinitesimal deformation functor of $C'_x$ is
  closed in the infinitesimal deformation functor of the new part of the
  relative jacobian of $C''_x/C'_x$, a generalised Prym variety.
\end{itemize}
  
A close inspection of \cite{Tama2} shows that we may moreover assume first,
that $C'_x \to C_x$ is $\ell_1$-primary for some arbitrary but sufficiently
large prime $\ell_1$, and then secondly the degree of $C''_x/C'_x$ can be
chosen to be a prime $\ell_2$, which is again arbitrary but sufficiently
large. It is exactly here that our approach uses the outer pro-$\LL$
representations for all pairs $\LL$ of two sufficiently large prime numbers
$\ell_1,\ell_2$.
  
By what was said above, the \'etale covers $C''_x \to C'_x \to C_x$ are
isomorphic to the fibres above $x$ of \'etale covers $C'' \to C' \to C$. The
cover $C''/C'$ remains prime cyclic Galois.  It follows that the new part of
the relative jacobian of the cyclic cover $C''/C'$ is a projective family of
ordinary abelian varieties, and hence constant by a theorem of Raynaud, Szpiro
and Moret-Bailly, see \cite{MBpinceauVarAb} XI Thm 5.1.  Consequently, by
property (iii) `new Torelli' above, the deformation $C'$ of $C'_x$ is a
constant family. A priori it is constant in the sense that all fibres are
mutually isomorphic, but as $C'$ carries a level structure this implies that
$C'/F$ is constant.  The original curve $C/F$ is a quotient of the Galois
closure of $C'/C$ and thus also constant, which proves Theorem \ref{tg}.
\end{absatz}
\end{pro}

An application of Theorem \ref{tc} is the following immediate corollary.

\begin{cor} \label{ta}
  Let $X$ be a reduced, connected and simply connected variety, i.e.
  $\pi_1X=1$, over an algebraically closed field $k$. Then any smooth proper
  curve $C$ over $X$ of genus $g$ with $n$ disjoint sections, where $2-2g-n$
  is negative, is a constant family of such curves, i.e, comes by base change
  from an object in $\mgn\big(\Spec(k)\big)$.
\end{cor}

\begin{absatz}[Remark.]
  In particular, Corollary \ref{ta} claims that there are no $\PP^1$'s on
  $\mgn$ for $2-2g-n < 0$. Note that the corresponding result for polarized
  abelian varieties fails. Oort has constructed a family of abelian varieties
  over $\PP^2_k$ with image of dimension $2$ in the moduli space of polarized
  abelian varieties, \cite{Oort} Rmk 2.6. Moreover, the coarse moduli space
  $M_g$ of proper smooth curves behaves also different from the fine moduli
  space $\mg$ that we use in this note.  Again Oort proves the existence of a
  complete rational curve in $M_g \otimes k$ for $k$ of positive
  characteristic and suitable $g$, \cite{Oort} Thm 3.1.
\end{absatz}


\section{Concluding remarks}

\begin{absatz}[Remarks.]
  (a) Let $S$ be a connected local complete intersection scheme and $U$ a
  dense open in $S$ such that $S-U$ has codimension at least $3$. In
  \cite{SGA2} Exp.\ X Thm 3.4 Grothendieck proves that $\pi_1 U \to \pi_1 S$
  is an isomorphism. Hence, if $S$ is normal, by Theorem \ref{t} the natural
  functor $\mgn(S)\to \mgn(U)$ is an equivalence of categories.
  
  (b) In order to yield a more natural anabelian proof of Theorem \ref{tc}, it
  would be nice if the natural map $\Ho^2(\pi_1\mgn,\Q_\ell) \to
  \Ho^2(\mgn,\Q_\ell)$ is an isomorphism. Note that $\Ho^2(\mgn,\Q_\ell) =
  \Q_\ell \cdot {\rm c}_1(\lambda)$. Then Theorem \ref{tc} allows a quick
  proof exploiting the $\ell$-adic first Chern class of the Hodge bundle
  $\lambda$, see \cite{stix2}.  In particular, this would give a proof of
  Theorem \ref{to} unconditionally on $\ell \in \OO_S^\ast$ for all
  characteristics.
  
  The isomorphism $\Ho^2(\pi_1\mgn,\Q_\ell) \cong \Ho^2(\mgn,\Q_\ell)$ has
  been announced by Boggi.
\end{absatz}

\begin{thm}[extended monodromy criterion] \label{ext}
  Let $2-2g-n$ be negative.  Let $f :U \to S$ be a map of finte type that
  factorizes as $f= \ov{f} \circ j$ where $j: U \to X$ is a dense open
  immersion into a normal scheme $X$, and $\ov{f}: X \to S$ is a proper map
  that satisfies $\ov{f}_\ast \OO_X = \OO_S$. Furthermore we assume that $S$
  is excellent. Then the following holds.
\begin{itemize}
\item[(1)] The pull back functor $f^\ast: \mgn(S) \to \mgn(U)$ is fully
  faithful.
\item[(2)] A $U$-curve $C \in \mgn(U)$ comes by base change from an $S$-curve
  in $\mgn(S)$ if and only if for all $N \in \N$ and geometric points $\ov{u}
  \in U \otimes \Z[\bru{0.5pt}{1}{N}]$ the following commutative diagram
  exists, where the solid arrows are induced by the natural maps $U \to \mgn$
  and $f: U \to S$:
  $$
  \xymatrix{{\pi_1(U \otimes \Z[\bru{0.5pt}{1}{N}],\ov{u})} \ar[rr] \ar[dr]
    &&
    {\pi_1(\mgn \otimes \Z[\bru{0.5pt}{1}{N}],C_{\ov{u}})} \\
    & {\pi_1(S \otimes \Z[\bru{0.5pt}{1}{N}],\ov{u})} \ar@{.>}[ur]& }$$
\end{itemize}
\end{thm}

\begin{pro}
  By applying Theorem \ref{t} to $U \subset X$ we may assume that $f$ is
  proper. Using again the finiteness of the isom-scheme of hyperbolic curves
  (resp. abelian level structures plus descent), (1) and (2) come down to the
  property of $\ov{f} : X \to S$ that a map of schemes $h : X \to T$
  factorizes through $\ov{f}$ if and only if it is constant along the fibres
  of $\ov{f}$ as a map of sets. In particular, we use again Theorem \ref{tc}
  for claim (2).
\end{pro}   

The obvious formulation of the analogous theorem in terms of outer pro-$\LL$
representations is left to the reader.


\end{document}